# PERFECT SPHERE PACKING IN THE BOOLEAN SPACE

T. SOGHOMONYAN[a,*], Zh. MARGARYAN[b]

[a] Faculty of Informatics and Applied Mathematics, Yerevan State University, 1 Alek Manukyan St., Yerevan 0025, Armenia
[b] Department of Discrete Mathematics and Theoretical Informatics, Faculty of Informatics and Applied Mathematics, Yerevan State University, 1 Alek Manukyan St., Yerevan 0025, Armenia

***Abstract.*** This article investigates the problem of perfect sphere packing of a Boolean space. Perfect sphere packing is a fundamental and complex problem with profound consequences for discrete mathematics, coding theory and cryptography. The classical solution to the perfect sphere packing problem was provided by Hamming via his well-known perfect codes. However, it is feasible only for spaces with specific dimensions. To address this limitation, this paper introduces new distance classes for which exists a perfect sphere packing for a Boolean space of arbitrary dimension. In this work, we explore the conditions for perfect Hamming packing based on a new distance and show that they hold for a space of arbitrary dimension. Furthermore, we construct an algorithm that provides perfect sphere packing for a Boolean space of arbitrary dimension.

## 1. Introduction

One of the problems in discrete mathematics, information and coding theory, is the perfect sphere packing of the Boolean space. In classical Hamming codes, such perfect sphere packing of the space is possible only for codes of strictly limited length, a limitation thoroughly investigated in the context of the nonexistence of perfect codes over finite fields [25]. From a discrete mathematics perspective, the essence of $q = 2$ block codes lies in partitioning the Boolean space into disjoint spheres or partitions [17]. If these spheres cover the entire space, the respective codes are considered perfect [4, 13, 15, 20]. The foundation of block codes was introduced by Hamming with his perfect codes [14], which remain widely used to this day in practical systems like semiconductor memory applications [7]. Over the decades, various alternative metrics and error models have been proposed to generalize or modify the classical Hamming space. These include the Lee metric and its associated perfect codes [3, 8, 12, 16, 18], which are closely related to multidimensional geometric coverings [23] and linear partitions [2]. Additionally, codes for asymmetric and additive communication channels have been extensively studied [1, 9, 11, 25], alongside other structured distances such as poset metrics [6] and combinatorial vector space problems [21, 22]. However, achieving a perfect sphere packing for an arbitrary dimension remains an open challenge [10].

The purpose of this article is to achieve perfect sphere packing for all dimensions in the Boolean space, thereby overcoming the limitations of the classical Hamming method. This is done by defining a new generalized Hamming distance, which provides a more respective mathematical generalization than existing models that rely on specific bit-shift operations [19].

* Corresponding author. Email: tigran.soghomonyan2@edu.ysu.am

## 2. HAMMING'S PERFECT SPHERE PACKING

Let $B^n = \{\boldsymbol{\alpha} |\ \boldsymbol{\alpha} = (\alpha_1, \alpha_2, \cdots, \alpha_n), \alpha_i \in \{0,1\}, i = \overline{1,n}\}$ be a space over the field $\mathbb{Z}_2$. For the distance between any arbitrary vectors $\boldsymbol{\alpha}, \boldsymbol{\beta} \in B^n$ in this space, we will use the classical Hamming distance. The Hamming distance is defined as the number of positions at which these vectors differ from each other:

$$d(\boldsymbol{\alpha}, \boldsymbol{\beta}) = \sum_{i=1}^{n} |\alpha_i - \beta_i|.$$

It satisfies all the metric axioms. Using the Hamming distance, we now give the definitions of a sphere and a ball in the space $B^n$.

**Definition 2.1** (Sphere)**.** A sphere of radius $t$ centered at $\boldsymbol{\alpha}$ in the space $B^n$ is defined as the following set of vectors:

$$\mathcal{S}_t(\boldsymbol{\alpha}) = \{\boldsymbol{\beta} \in B^n \mid d(\boldsymbol{\alpha}, \boldsymbol{\beta}) = t\}.$$

**Definition 2.2** (Ball)**.** A ball of radius $t$ centered at $\boldsymbol{\alpha}$ in the space $B^n$ is defined as the following set of vectors:

$$\mathcal{B}_t(\boldsymbol{\alpha}) = \{\boldsymbol{\beta} \in B^n \mid d(\boldsymbol{\alpha}, \boldsymbol{\beta}) \leq t\}.$$

**Definition 2.3** (Minimal distance)**.** Let $L \subseteq B^n$ is an arbitrary subset. The minimal distance of $L$ is defined as the smallest distance between any two distinct vectors in this subset:

$$d = \min_{\substack{\boldsymbol{\alpha}, \boldsymbol{\beta} \in L, \\ \boldsymbol{\alpha} \neq \boldsymbol{\beta}}} \{d(\boldsymbol{\alpha}, \boldsymbol{\beta})\}.$$

**Definition 2.4** (Hamming weight)**.** The weight of a vector $\boldsymbol{\alpha} \in B^n$ belonging to the Boolean space is defined as the number of its non-zero components:

$$w(\boldsymbol{\alpha}) = \sum_{i=1}^{n} \alpha_i.$$

In other words, the weight of a vector can be represented as the distance between it and the zero vector: $w(\boldsymbol{\alpha}) = d(\boldsymbol{\alpha}, \mathbf{0})$. It is well known from coding theory that the minimum distance of a space is equal to the minimum weight of its non-zero vectors:

$$d = \min_{\substack{\boldsymbol{\gamma} \in L, \\ \boldsymbol{\gamma} \neq \mathbf{0}}} \{w(\boldsymbol{\gamma})\}.$$

Let $C \subseteq B^n$ be an any subspace. For the space $B^n$ to be partitioned into balls of radius $t$, where the centers of the balls are vectors of the subspace $C$, it is necessary that the minimum distance of the subspace $C$ satisfies the $d = 2t + 1$ condition. This condition ensures that the balls formed around the vectors of the subspace $C$ do not intersect pairwise. In addition, it is necessary that the sum of the cardinalities of these balls is equal to the cardinality of the space $B^n$. This sum cannot exceed the cardinality of the space $B^n$, which is known as the Hamming upper packing bound [4, 14, 15, 20]. However, recent works have proposed further refinements to this bound for binary codes [21]. The codes for which the upper packing bound becomes an equality are called perfect codes [15, 17]. The most famous solution to this problem, in the case of radius $t = 1$, is the Hamming perfect codes [14].

It is known from the theory of linear codes that a linear code $C$ can be given by its parity-check matrix $\boldsymbol{H}$ [15].

**Theorem 2.1.** Let $\boldsymbol{H}_{(m\times n)}$ be the parity-check matrix of the linear code $C$. For the minimum distance of the linear code $C$ to be equal to $d$, it is necessary and sufficient that the columns of the parity-check matrix $\boldsymbol{H}$ satisfy the following two conditions:

- $\forall d-1$ columns are linearly independent.
- $\exists d$ columns that are linearly dependent.

The proof of the theorem can be found in [15]. According to Theorem 2.1, for the minimum distance of the code to be exactly $d=3$, it is necessary and sufficient that any two columns of the parity-check matrix $\boldsymbol{H}$ are linearly independent and there exist three columns that are linearly dependent. In the field $\mathbb{Z}_2$, the linear independence of two columns means that all columns of the matrix $\boldsymbol{H}$ must be non-zero and pairwise distinct. Since the columns of the matrix have dimension $m$, the maximum number of non-zero and distinct $m$-dimensional vectors is equal to $2^m-1$. All non-zero vectors are selected as the columns of the parity-check matrix $\boldsymbol{H}$ of the Hamming perfect codes. Such a construction ensures the conditions of Theorem 2.1. Since the zero column is absent in the matrix and all columns are pairwise distinct, therefore any two of them are linearly independent. In addition, the matrix contains all possible non-zero column vectors of dimension $m$, therefore the sum of any two distinct columns will result in another non-zero column vector, which is also a column of the matrix. As a result, there always exist three columns that are linearly dependent. On the other hand, the sum of the cardinalities of the balls is equal to the number of vectors in the space $B^n$ which shows that the Hamming upper packing bound is satisfied with equality.

In the case of the classical Hamming distance, a perfect sphere packing of the space $B^n$ with balls of radius one is possible if and only if the dimension of the space is $n=2^m-1$. For any other dimensions, it is impossible to obtain a perfect sphere packing with the classical metric.

## 3. GENERALIZED HAMMING DISTANCE

The partition method is based on the generalized distance, based on the Hamming distance. Consider the set of indices $N=\{1,2,\cdots,n\}$ and fix the subsets $I_1,I_2,\cdots,I_k\subseteq N$ for which it holds that:

$$\bigcup_{j=1}^{k} I_j = N.$$

**Definition 3.1** (Partial distance)**.** For a given subset of indices $I_j$, $\boldsymbol{\alpha}$ and $\boldsymbol{\beta}$ is defined as:

$$d_j(\boldsymbol{\alpha},\boldsymbol{\beta}) = \sum_{i\in I_j} |\alpha_i-\beta_i|.$$

**Definition 3.2** (Generalized distance)**.** For fixed subsets of indices $I_1,I_2,\cdots,I_k$ and arbitrary vectors $\boldsymbol{\alpha}$ and $\boldsymbol{\beta}$, the generalized distance is defined as the maximum value among the partial distances over the individual subsets:

$$D(\boldsymbol{\alpha},\boldsymbol{\beta}) = \max\{d_1(\boldsymbol{\alpha},\boldsymbol{\beta}), d_2(\boldsymbol{\alpha},\boldsymbol{\beta}), \cdots, d_k(\boldsymbol{\alpha},\boldsymbol{\beta})\}.$$

**Proposition 3.1:** The generalized distance satisfies the four axioms of a metric, for arbitrary vectors $\boldsymbol{\alpha},\boldsymbol{\beta},\boldsymbol{\omega} \in B^n$ the following statements hold:

1. **Non-negativity:** $D(\boldsymbol{\alpha},\boldsymbol{\beta}) \geq 0$.
2. **Identity:** $D(\boldsymbol{\alpha},\boldsymbol{\beta}) = 0 \Leftrightarrow \boldsymbol{\alpha} = \boldsymbol{\beta}$.
3. **Identity:** $D(\boldsymbol{\alpha},\boldsymbol{\beta}) = D(\boldsymbol{\beta},\boldsymbol{\alpha})$.
4. **Triangle inequality:** $D(\boldsymbol{\alpha},\boldsymbol{\beta}) \leq D(\boldsymbol{\alpha},\boldsymbol{\omega}) + D(\boldsymbol{\omega},\boldsymbol{\beta})$.

**Proof.**

1. **Non-negativity:** By definition, each $d_j(\boldsymbol{\alpha},\boldsymbol{\beta})$ represents a sum of absolute values. Therefore, it is a non-negative quantity: $d_j(\boldsymbol{\alpha},\boldsymbol{\beta}) \geq 0$ for all indices $1 \leq j \leq k$. Therefore, their maximum value is also non-negative:
$$D(\boldsymbol{\alpha},\boldsymbol{\beta}) = \max_{1\leq j\leq k} d_j(\boldsymbol{\alpha},\boldsymbol{\beta}) \geq 0.$$
2. **Identity:** Let us assume that $D(\boldsymbol{\alpha},\boldsymbol{\beta}) = 0$ and show that $\boldsymbol{\alpha} = \boldsymbol{\beta}$. If
$$\max_{1\leq j\leq k} d_j(\boldsymbol{\alpha},\boldsymbol{\beta}) = 0,$$
then taking into account the non-negativity, it follows that for each $j$, $d_j(\boldsymbol{\alpha},\boldsymbol{\beta}) = 0$. Since the subsets $I_1, I_2, \cdots, I_k$ cover the entire set of indices, then for any index we obtain $|\alpha_i - \beta_i| = 0$, from which it is obtained that $\boldsymbol{\alpha} = \boldsymbol{\beta}$. Let us show the other direction. If $\boldsymbol{\alpha} = \boldsymbol{\beta}$, then for an arbitrary index $|\alpha_i - \beta_i| = 0$. Therefore, for any subset $I_j$, $d_j(\boldsymbol{\alpha},\boldsymbol{\beta}) = 0$. If it is true for each of them, then it is also true for their maximum value: $D(\boldsymbol{\alpha},\boldsymbol{\beta}) = 0$.
3. **Symmetry:** Since the distance is based on the absolute difference and $|\alpha_i - \beta_i| = |\beta_i - \alpha_i|$, it is obvious that $d_j(\boldsymbol{\alpha},\boldsymbol{\beta}) = d_j(\boldsymbol{\beta},\boldsymbol{\alpha})$. From which it follows that $D(\boldsymbol{\alpha},\boldsymbol{\beta}) = D(\boldsymbol{\beta},\boldsymbol{\alpha})$.
4. **Triangle inequality:** Let us consider an arbitrary subset of indices $I_j$ $(1 \leq j \leq k)$. For an arbitrary index $i \in I_j$ of that set, it is true that $|\alpha_i - \beta_i| \leq |\alpha_i - \omega_i| + |\omega_i - \beta_i|$. Summing this inequality for all indices $i \in I_j$, we will obtain:
$$\sum_{i\in I_j}|\alpha_i - \beta_i| \leq \sum_{i\in I_j}(|\alpha_i - \omega_i| + |\omega_i - \beta_i|) = \sum_{i\in I_j}|\alpha_i - \omega_i| + \sum_{i\in I_j}|\omega_i - \beta_i|.$$
According to our notations, we obtain $d_j(\boldsymbol{\alpha},\boldsymbol{\beta}) \leq d_j(\boldsymbol{\alpha},\boldsymbol{\omega}) + d_j(\boldsymbol{\omega},\boldsymbol{\beta})$. Since the inequality is true for an arbitrary $I_j$, then it is also true for the maximum value:

$$\max_{1\leq j\leq k} d_j(\boldsymbol{\alpha},\boldsymbol{\beta}) \leq \max_{1\leq j\leq k}\{d_j(\boldsymbol{\alpha},\boldsymbol{\omega}) + d_j(\boldsymbol{\omega},\boldsymbol{\beta})\} \leq \max_{1\leq j\leq k} d_j(\boldsymbol{\alpha},\boldsymbol{\omega}) + \max_{1\leq j\leq k} d_j(\boldsymbol{\omega},\boldsymbol{\beta}) \blacksquare$$

**Remark 3.1:** Note that the generalized distance in a special case reduces to the classical Hamming distance. In the case when at least one of the selected subsets of indices $I_1, I_2, \cdots, I_k$ coincides with the entire set $N$ (that is, there exists such a $j \in \{1, \cdots, k\}$, that $I_j = N$), then the partial distances calculated on all other subsets cannot exceed $d_j$. As a result, the generalized distance coincides with the Hamming distance.

**Example 3.1:** Let $N = \{1, 2, 3, 4, 5\}$, and the covering is given:

$$I_1 = \{2, 3, 4\}, I_2 = \{1, 2\}, I_3 = \{5\}.$$

Let us consider two vectors: $\boldsymbol{\alpha} = (0, 1, 0, 1, 0)$ and $\boldsymbol{\beta} = (1, 0, 1, 1, 0)$. Let us calculate the distances in the subsets $I_j$:

$$d_1(\boldsymbol{\alpha}, \boldsymbol{\beta}) = |\alpha_2 - \beta_2| + |\alpha_3 - \beta_3| + |\alpha_4 - \beta_4| = 1 + 1 + 0 = 2,$$

$$d_2(\boldsymbol{\alpha}, \boldsymbol{\beta}) = |\alpha_1 - \beta_1| + |\alpha_2 - \beta_2| = 1 + 1 = 2,$$

$$d_3(\boldsymbol{\alpha}, \boldsymbol{\beta}) = |\alpha_5 - \beta_5| = 0,$$

$$D(\boldsymbol{\alpha}, \boldsymbol{\beta}) = \max\{d_1(\boldsymbol{\alpha}, \boldsymbol{\beta}),\ d_2(\boldsymbol{\alpha}, \boldsymbol{\beta}), d_3(\boldsymbol{\alpha}, \boldsymbol{\beta})\} = 2.$$

Note that the classical Hamming distance for these vectors is $d(\boldsymbol{\alpha}, \boldsymbol{\beta}) = 3$. This example shows that the generalized distance can be less than the Hamming distance.

The definitions of a sphere and a ball do not undergo a change, simply instead of the Hamming distance, the new generalized distance $D$ is applied. The sphere and the ball with an arbitrary center $\boldsymbol{\alpha}$ and radius $r$ in the space $B^n$, according to the generalized distance $D$ will have the following form:

$$\mathcal{S}_r(\boldsymbol{\alpha}) = \{\boldsymbol{\beta} \in B^n \mid D(\boldsymbol{\alpha}, \boldsymbol{\beta}) = r\}, \mathcal{B}_r(\boldsymbol{\alpha}) = \{\boldsymbol{\beta} \in B^n \mid D(\boldsymbol{\alpha}, \boldsymbol{\beta}) \leq r\}.$$

The principle of calculating the cardinalities of a sphere and a ball with the new metric differs significantly. Since the ball represents the union of all spheres located up to a distance $r$, its cardinality can be calculated as the sum of the cardinalities of the corresponding spheres:

$$|\mathcal{B}_r(\boldsymbol{\alpha})| = \sum_{i=0}^{r} |\mathcal{S}_i(\boldsymbol{\alpha})|.$$

Thus, the study of the cardinality of a ball reduces to the cardinality of a sphere with a fixed radius.

The formula for the cardinality of a sphere is closely connected to the subsets of indices. Before the derivation of the formula for the cardinality of a sphere, it is necessary to perform some preparatory work.

**Definition 3.3** (Partial weight)**:** Let $I \subseteq N$ be an arbitrary subset of indices, and $\boldsymbol{\alpha} \in B^n$ be an arbitrary vector. The weight of the vector $\boldsymbol{\alpha}$ in the subset of indices $I$ will be defined as follows:

$$w_I(\boldsymbol{\alpha}) = \sum_{i \in I} \alpha_i.$$

**Definition 3.4** (Generalized weight)**:** Let $I_1, I_2, \cdots, I_k \subseteq N$ be arbitrary subsets of indices, and $\boldsymbol{\alpha} \in B^n$ be an arbitrary vector. The weight of the vector $\boldsymbol{\alpha}$ in the metric $D$ will be defined as follows:

$$w_D(\boldsymbol{\alpha}) = \max\{w_{I_1}(\boldsymbol{\alpha}), w_{I_2}(\boldsymbol{\alpha}), \cdots, w_{I_k}(\boldsymbol{\alpha})\}.$$

*Further analyses will be carried out for fixed subsets of indices: $I_1$ and $I_2$.*

Let us partition the set $N$ into three disjoint parts, and introduce the following notations:

- $A \equiv N \setminus I_2$, the number of elements of which is $n_A \equiv |A|$,
- $B \equiv N \setminus I_1$, the number of elements of which is $n_B \equiv |B|$,
- $C \equiv I_1 \cap I_2$, the number of elements of which is $n_C \equiv |C|$.

It is obvious that for the vector $\forall \boldsymbol{\alpha} \in B^n$, the following bounds hold:

$$0 \leq w_A(\boldsymbol{\alpha}) \leq n_A,\ 0 \leq w_B(\boldsymbol{\alpha}) \leq n_B,\ 0 \leq w_C(\boldsymbol{\alpha}) \leq n_C.$$

By applying these notations, the weights determined by the subsets $I_1$ and $I_2$ will be represented as the sum of partial weights:

$$w_{I_1}(\boldsymbol{\alpha}) = w_A(\boldsymbol{\alpha}) + w_C(\boldsymbol{\alpha}),\ w_{I_2}(\boldsymbol{\alpha}) = w_B(\boldsymbol{\alpha}) + w_C(\boldsymbol{\alpha}).$$

And the distance with the new metric will take the following form:

$$D(\boldsymbol{\alpha}, \boldsymbol{\beta}) = \max\{w_{I_1}(\boldsymbol{\alpha} \oplus \boldsymbol{\beta}), w_{I_2}(\boldsymbol{\alpha} \oplus \boldsymbol{\beta})\}.$$

**Proposition 3.2** (Cardinality of a sphere)**:** In the space $B^n$, the cardinality of a sphere with an arbitrary center $\boldsymbol{\alpha}$ and radius $r$, in the $D$ metric, is determined by the following formula:

$$|S_r(\boldsymbol{\alpha})| = S_1 + S_2 - S_{12}, \text{ where:}$$

$$S_1 = \sum_{w_C=0}^{\min\{n_C,r\}} \left[ \binom{n_C}{w_C} \binom{n_A}{r - w_C} \sum_{w_B=0}^{\min\{n_B, r-w_C\}} \binom{n_B}{w_B} \right],$$

$$S_2 = \sum_{w_C=0}^{\min\{n_C,r\}} \left[ \binom{n_C}{w_C} \binom{n_B}{r - w_C} \sum_{w_A=0}^{\min\{n_A, r-w_C\}} \binom{n_A}{w_A} \right],$$

$$S_{12} = \sum_{w_C=0}^{\min\{n_C,r\}} \binom{n_C}{w_C} \binom{n_A}{r - w_C} \binom{n_B}{r - w_C}.$$

**Proof.** Let us consider the condition $\max\{w_{I_1}(\boldsymbol{\alpha} \oplus \boldsymbol{\beta}), w_{I_2}(\boldsymbol{\alpha} \oplus \boldsymbol{\beta})\} = r$. It can be represented as the union of the following two conditions:

1. $\left(w_{I_1}(\boldsymbol{\alpha} \oplus \boldsymbol{\beta}) = r \textit{ and } w_{I_2}(\boldsymbol{\alpha} \oplus \boldsymbol{\beta}) \leq r\right)$,
2. $\left(w_{I_2}(\boldsymbol{\alpha} \oplus \boldsymbol{\beta}) = r \textit{ and } w_{I_1}(\boldsymbol{\alpha} \oplus \boldsymbol{\beta}) \leq r\right)$.

Therefore, to calculate the total number, it is necessary to consider these cases separately, and according to the inclusion-exclusion principle, subtract the case of simultaneous intersection $\left(w_{I_1} = r \text{ and } w_{I_2} = r\right)$ from the sum.

**Case one.** Let us consider the vectors for which the first condition is satisfied: $w_{I_1} = r$ and $w_{I_2} \leq r$. Let us denote the number of those vectors by $S_1$.

Since according to our notations $w_{I_1} = w_A + w_C$ and $w_{I_2} = w_B + w_C$, let us make the selections sequentially over the subsets $C$, $A$ and $B$.

1. **Subset $C$:** First, let us choose the number of weights $w_C$ in the $n_C$ positions of the subset $C$. Since the condition $w_C \leq n_C$ must be ensured, the number of possible variants is $\binom{n_C}{w_C}$, the value of $w_C$ varies from 0 to $\min\{n_C, r\}$:
2. **Subset $A$:** In order to ensure the condition $w_{I_1} = r$, there must be exactly $r - w_C$ ones in the $n_A$ positions of the subset $A$. The number of choices is $\binom{n_A}{r-w_C}$. This factor is meaningful and non-zero only in the case when $r - w_C \leq n_C$.
3. **Subset $B$:** In order to ensure the condition $w_{I_2} \leq r$, the number of ones in the subset $B$ can vary from 0 to $r - w_C$. Therefore, the total number of choices will be obtained in the form of a sum:

$$\sum_{w_B=0}^{\min\{n_B, r-w_C\}} \binom{n_B}{w_B}.$$

Combining these three steps, we will obtain the final formula for $S_1$:

$$S_1 = \sum_{w_C=0}^{\min\{n_C, r\}} \left[ \binom{n_C}{w_C} \binom{n_A}{r - w_C} \sum_{w_B=0}^{\min\{n_B, r-w_C\}} \binom{n_B}{w_B} \right].$$

**Case two.** Consider the vectors for which the second condition is satisfied: $w_{I_2} = r$ and $w_{I_1} \leq r$. It is symmetric to the first one; we will denote the number of vectors satisfying that condition by $S_2$. By similar reasoning, we will obtain:

$$S_2 = \sum_{w_C=0}^{\min\{n_C, r\}} \left[ \binom{n_C}{w_C} \binom{n_B}{r - w_C} \sum_{w_A=0}^{\min\{n_A, r-w_C\}} \binom{n_A}{w_A} \right].$$

**Case three.** Consider the case of simultaneous occurrence of the two conditions (the intersection): $w_{I_1} = r$ and $w_{I_2} = r$. We will denote the number of vectors satisfying this condition by $S_{12}$. Since both conditions require a weight of exactly $r$, then after choosing $w_C$ positions from the subset $C$, exactly $r - w_C$ positions must be chosen from the remaining subsets $A$ and $B$:

$$S_{12} = \sum\nolimits_{w_C} \binom{n_C}{w_C} \binom{n_A}{r - w_C} \binom{n_B}{r - w_C}.$$

Thus, according to the inclusion-exclusion principle, the cardinality of the sphere defined under the new metric will be calculated by the following formula:

$$\max\{0, r - n_A, r - n_B\} \leq w_C \leq \min\{n_C, r\}.$$

Thus, by the inclusion-exclusion principle, the cardinality of the sphere defined under the new metric is calculated by the following formula:

$$|\mathcal{S}_r(\boldsymbol{\alpha})| = S_1 + S_2 - S_{12}. \blacksquare$$

**Example 3.2:** For simplicity, consider the following example:

$N = \{1,2,3,4,5\}, I_1 = \{1,2,3,4\}, I_2 = \{3,4,5\}$

$A = \{1,2\}, n_A = 2$

$B = \{5\}, n_B = 1$

$C = \{3,4\}, n_C = 2$

Calculate $|\mathcal{S}_r(\boldsymbol{\alpha})|$, where $\boldsymbol{\alpha} \in B^5$ and $r = 2$:

$$S_1 = \sum_{w_C=0}^{\min\{2,2\}} \left[ \binom{2}{w_C} \binom{2}{2-w_C} \sum_{w_B=0}^{\min\{1,2-w_C\}} \binom{1}{w_B} \right]$$

$w_C = 0$: $\quad \binom{2}{0}\binom{2}{2} \sum_{w_B=0}^{\min\{1,2\}} \binom{1}{w_B} = 1 \cdot 1 \cdot \left(\binom{1}{0} + \binom{1}{1}\right) = 2$

$w_C = 1$: $\quad \binom{2}{1}\binom{2}{1} \sum_{w_B=0}^{\min\{1,1\}} \binom{1}{w_B} = 2 \cdot 2 \cdot \left(\binom{1}{0} + \binom{1}{1}\right) = 8$

$w_C = 2$: $\quad \binom{2}{2}\binom{2}{2} \sum_{w_B=0}^{\min\{1,0\}} \binom{1}{w_B} = 1 \cdot 1 \cdot \binom{1}{0} = 1$

$$S_1 = 8 + 2 + 1 = 11$$

$$S_2 = \sum_{w_C=0}^{\min\{2,2\}} \left[ \binom{2}{w_C} \binom{1}{2-w_C} \sum_{w_A=0}^{\min\{2,2-w_C\}} \binom{2}{w_A} \right]$$

$w_C = 0$: $\quad \binom{2}{0}\binom{1}{2} \sum_{w_A=0}^{\min\{2,2\}} \binom{2}{w_A} = 1 \cdot 0 \cdot \left(\binom{2}{0} + \binom{2}{1} + \binom{2}{2}\right) = 0$

$w_C = 1$: $\quad \binom{2}{1}\binom{1}{1} \sum_{w_A=0}^{\min\{2,1\}} \binom{2}{w_A} = 2 \cdot 1 \cdot \left(\binom{2}{0} + \binom{2}{1}\right) = 6$

$w_C = 2$: $\quad \binom{2}{2}\binom{1}{0} \sum_{w_A=0}^{\min\{2,0\}} \binom{2}{w_A} = 1 \cdot 1 \cdot \binom{2}{0} = 1$

$$S_2 = 0 + 6 + 1 = 7$$

$$S_{12} = \sum_{w_C=\max\{0,2-2,2-1\}}^{\min\{2,2\}} \binom{2}{w_C} \binom{2}{2-w_C} \binom{1}{2-w_C}$$

$w_C = 1$: $\quad \binom{2}{1}\binom{2}{1}\binom{1}{1} = 2 \cdot 2 \cdot 1 = 4$

$w_C = 2$: $\quad \binom{2}{2}\binom{2}{0}\binom{1}{0} = 1 \cdot 1 \cdot 1 = 1$

$$S_{12} = 4 + 1 = 5$$

$$|\mathcal{S}_2(\boldsymbol{\alpha})| = S_1 + S_2 - S_{12} = 11 + 7 - 5 = 13.$$

In this article, we will mainly deal with balls of radius one. Since the calculation of the cardinality of balls is based on the cardinality of spheres, it is worth noting a simpler method for calculating the cardinality of a sphere of radius one. According to the new metric, vectors at a distance of one (radius one) can occur in two cases:

1. When the vectors differ from each other in exactly one position (the number of such possible variants is $n$).
2. When the vectors differ from each other in exactly one position (the number of such possible variants is $A$, and the other to $B$ (the number of choices for such pairs is $n_A \cdot n_B$).

Therefore, the cardinality of a sphere of radius one is calculated by the following simple formula:

$$|\mathcal{S}_1(\boldsymbol{\alpha})| = n + n_A \cdot n_B.$$

In order for it to be possible to partition the space into balls, it is necessary that the cardinalities of the balls be powers of two. Since the cardinality of a ball of radius one is the sum of the cardinalities of the spheres of radius 0 and 1, it will take the following form:

$$|\mathcal{B}_1(\boldsymbol{\alpha})| = |\mathcal{S}_0(\boldsymbol{\alpha})| + |\mathcal{S}_1(\boldsymbol{\alpha})| = 1 + n + n_A \cdot n_B.$$

Consequently, for the partition, there must exist integers $n_A$, $n_B$ and $m$ such that the following holds:

$$n + n_A \cdot n_B \ + \ 1 = 2^m \Rightarrow n_A \cdot n_B = 2^m \ - n - 1. \quad (1)$$

**Theorem 3.1** (Existence of partition parameters)**:** Under the new $D$ metric, for an arbitrary Boolean space of dimension $n \geq 2$, there exist parameters $n_A$ and $n_B$ such that the cardinality of balls of radius one satisfies the Hamming perfect packing condition.

**Proof.** For the Hamming perfect packing condition to be satisfied, there must exist a natural number $m \geq 2$ such that equation (1) holds, from which we obtain: $n_A \cdot n_B = 2^m \ - n - 1$. For the partition to be correct, the selected values of $n_A$ and $n_B$ must satisfy the condition $n_A + n_B \leq n$ (since they are the cardinalities of disjoint subsets of a set of length $n$).

To ensure this condition, it is sufficient to choose the value $m$ such that the resulting integer $2^m - n - 1$ can be factored into two multipliers whose sum does not exceed $n$. Regardless of whether the resulting number is prime or composite, it can always be represented by integer multipliers (in the case of a prime number, one of the multipliers will be 1). Let us show that it is always possible to find such a natural number $m$ for which both conditions will be satisfied.

For an arbitrary natural number $n \geq 1$, in the interval $(n, 2n]$, there always exists a power of two. That is, we can always choose a natural number $m$ for which the following inequality holds: $n + 1 \leq 2^m \leq 2n$. Selecting exactly this $m$, we can take $n_A = 1$ and $n_B = 2^m - n - 1$ as multipliers. For the sum of the multipliers, we will obtain:

$$n_A + n_B = 1 + 2^m - n - 1 = 2^m - n.$$

Since the inequality $2^m \leq 2n$ holds for our selected $m$, the required condition $n_A + n_B \leq n$ is always satisfied. ■

According to Theorem 3.1, here will always be a distance class for which the necessary condition for partition holds, but it does not fully guarantee a partition. For the partition of the space by balls of radius $t = 1$, it is also required that the minimal distance of the space, under the new $D$ metric, be $d_D = 3$ [13].

**Theorem 3.2** (Minimum distance)**:** Let $L \subset B^n$ be a linear subspace given by a parity-check matrix $\boldsymbol{H}$ and $I_1, I_2 \subseteq N$ be subsets of indices. For the minimum distance of the subspace $L$ under the $D$ metric to be exactly $d_D = 3$, it is necessary and sufficient that the columns of the parity-check matrix $\boldsymbol{H}$ satisfy the following conditions:

1. $\forall i \in N : \boldsymbol{h}_i \neq \boldsymbol{0}$ and $\forall i \neq j : \boldsymbol{h}_i \neq \boldsymbol{h}_j$,
2. $\exists i, j, k \in I_1$ or $\exists i, j, k \in I_2 : \boldsymbol{h}_i \oplus \boldsymbol{h}_j \oplus \boldsymbol{h}_k = \boldsymbol{0}$,
3. $\forall a \in A, \forall b \in B, \forall i \in N : \boldsymbol{h}_a \oplus \boldsymbol{h}_b \neq \boldsymbol{h}_i$,
4. $\forall a_1, a_2 \in A, \forall b_1, b_2 \in B : \boldsymbol{h}_{a_1} \oplus \boldsymbol{h}_{b_1} \neq \boldsymbol{h}_{a_2} \oplus \boldsymbol{h}_{b_2}$.

**Proof.** Since $L$ is a linear subspace, the minimum distance of the subspace under the $D$ metric is equal to the minimum generalized weight of a non-zero vector belonging to the subspace:

$$d_D = \min_{\substack{\boldsymbol{\alpha} \in L, \\ \boldsymbol{\alpha} \neq \boldsymbol{0}}} \{w_D(\boldsymbol{\alpha})\}.$$

**Sufficiency:** Assume conditions 1-4 hold, but there exists $\boldsymbol{\alpha} \neq \boldsymbol{0}$ such that $w_D(\boldsymbol{\alpha}) < 3$. Since $\boldsymbol{\alpha} \neq \boldsymbol{0}$, then $w_D(\boldsymbol{\alpha})$ can be equal to 1 or 2. Suppose $w_D(\boldsymbol{\alpha}) = 1$, by definition: $\max\{w_{I_1}(\boldsymbol{\alpha}), w_{I_2}(\boldsymbol{\alpha})\} = 1$. Such a vector $\boldsymbol{\alpha}$ can have three forms:

- $(w_{I_1} = 1, w_{I_2} = 0)$: $\boldsymbol{\alpha}$ has exactly one 1 in $A$, from which it follows that $\boldsymbol{H} \cdot \boldsymbol{\alpha}^T = \boldsymbol{h}_i = \boldsymbol{0}$, which contradicts the 1st condition.
- $(w_{I_1} = 0, w_{I_2} = 1)$: $\boldsymbol{\alpha}$ has exactly one 1 in $B$, this case is symmetric to the previous one and similarly contradicts the 1st condition.
- $(w_{I_1} = 1, w_{I_2} = 1)$: In this case, $\boldsymbol{\alpha}$ can have exactly one 1 in $C$ (at position $c$), from which it follows that $\boldsymbol{H} \cdot \boldsymbol{\alpha}^T = \boldsymbol{h}_c = \boldsymbol{0}$, which contradicts the 1st condition. The other possible variant is that the vector has a 1 in $A$ (at position $a$) and $B$ (at position $b$), from which it follows that $\boldsymbol{H} \cdot \boldsymbol{\alpha}^T = \boldsymbol{h}_a \oplus \boldsymbol{h}_b = \boldsymbol{0}$. This means $\boldsymbol{h}_a = \boldsymbol{h}_b$, which also contradicts the 1st condition.

Suppose $w_D(\boldsymbol{\alpha}) = 2$, from which five subcases arise:

- $(w_{I_1} = 2, w_{I_2} = 0)$: $\boldsymbol{\alpha}$ has two 1s only in $A$ (at positions $a_1$ and $a_2$), from which it follows that $\boldsymbol{H} \cdot \boldsymbol{\alpha}^T = \boldsymbol{h}_{a_1} \oplus \boldsymbol{h}_{a_2} = \boldsymbol{0}$, which contradicts the 1st condition.
- $(w_{I_1} = 0, w_{I_2} = 2)$: This case is symmetric to the previous one and contradicts the 1st condition.
- $(w_{I_1} = 2, w_{I_2} = 1)$: In this case, $\boldsymbol{\alpha}$ can have one 1 in $A$ (at position $a$) and one 1 in $C$ (at position $c$), from which it follows that $\boldsymbol{H} \cdot \boldsymbol{\alpha}^T = \boldsymbol{h}_a \oplus \boldsymbol{h}_c = \boldsymbol{0}$, contradicting the 1st condition. The other possible variant is that the vector has two 1s in $A$ (at positions $a_1$ and $a_2$) and one 1 in $B$ (at position $b$), from which it follows that $\boldsymbol{H} \cdot \boldsymbol{\alpha}^T = \boldsymbol{h}_{a_1} \oplus \boldsymbol{h}_{a_2} \oplus \boldsymbol{h}_b = \boldsymbol{0}$, which contradicts the 3rd condition.
- $(w_{I_1} = 1, w_{I_2} = 2)$: This case is symmetric to the previous one and contradicts the 1st or 3rd condition.
- $(w_{I_1} = 2, w_{I_2} = 2)$: In this case, $\boldsymbol{\alpha}$ can have two 1s $C$ (at positions $c_1$ and $c_2$), from which it follows that $\boldsymbol{H} \cdot \boldsymbol{\alpha}^T = \boldsymbol{h}_{c_1} \oplus \boldsymbol{h}_{c_2} = \boldsymbol{0}$, contradicting the 1st condition. The other possible variants are: a single 1 in all three subsets (contradicting the 3rd condition), or two 1s in both subsets $A$ and $B$ (at positions $a_1$, $a_2$, $b_1$ and $b_2$), from which it follows that $\boldsymbol{H} \cdot \boldsymbol{\alpha}^T = \boldsymbol{h}_{a_1} \oplus \boldsymbol{h}_{a_2} \oplus \boldsymbol{h}_{b_1} \oplus \boldsymbol{h}_{b_2} = \boldsymbol{0}$, which contradicts the 4th condition.

In all cases, we obtained a contradiction, therefore $d_D \geq 3$. Condition 2 ensures that there exists at least one vector $\boldsymbol{\alpha}$ for which $w_D(\boldsymbol{\alpha}) = 3$, excluding the case $d_D > 3$. Consequently, $d_D = 3$.

**Necessity :** Assume the contrary: $d_D = 3$ (i.e., for all vectors $\boldsymbol{\alpha} \in L \setminus \{\boldsymbol{0}\}$, $w_D(\boldsymbol{\alpha}) \geq 3$), but at least one of the four specified conditions does not hold. If any of conditions 1, 3, or 4 is violated, we will obtain a linear combination of the columns of matrix $\boldsymbol{H}$ that is equal to the zero-column vector, forming a codeword with $w_D < 3$. If the 1st condition is violated ($\boldsymbol{h}_i = \boldsymbol{0}$ or $\boldsymbol{h}_i = \boldsymbol{h}_j$), then there will exist a vector $\boldsymbol{\alpha}$ whose generalized weight will be less than 3 ($w_D(\boldsymbol{\alpha}) < 3$). If the 3rd condition is violated ($\boldsymbol{h}_a \oplus \boldsymbol{h}_b = \boldsymbol{h}_i$), then $\boldsymbol{h}_a \oplus \boldsymbol{h}_b \oplus \boldsymbol{h}_i = \boldsymbol{0}$. This means that there exists a vector $\boldsymbol{\alpha}$ for which $w_{I_1}(\boldsymbol{\alpha}) \leq 2$ or $w_{I_2}(\boldsymbol{\alpha}) \leq 2$. Therefore:

$$w_D(\boldsymbol{\alpha}) = \max\{w_{I_1}(\boldsymbol{\alpha}), w_{I_2}(\boldsymbol{\alpha})\} \leq 2.$$

If the 4th condition is violated ($\boldsymbol{h}_{a_1} \oplus \boldsymbol{h}_{b_1} = \boldsymbol{h}_{a_2} \oplus \boldsymbol{h}_{b_2}$), then the sum of these columns will be equal to the zero column vector. In this case as well, there will exist a vector $\boldsymbol{\alpha}$ for which $w_{I_1}(\boldsymbol{\alpha}) = 2$ and $w_{I_2}(\boldsymbol{\alpha}) = 2$. If the 2nd condition is violated, then there will be no vector $\boldsymbol{\alpha}$ found in the subspace for which $w_D(\boldsymbol{\alpha}) = 3$ would be obtained. Since satisfying the 1st, 3rd, and 4th conditions already prohibits the existence of code vectors with a generalized weight less than 3, violating the 2nd condition would cause the minimum distance to become strictly greater than 3. As a result, we obtain a contradiction. ■

## 4. Construction Of The Sphere Packing

In coding theory, it often happens that the theoretical existence of a matrix satisfying certain conditions is proven, but finding it is a computationally complex problem. For this reason, the theorem on the existence of the sphere packing is formulated below.

**Theorem 4.1** (Existence of the sphere packing)**:** An arbitrary Boolean space of dimension $n \geq 3$ can be partitioned into a sphere packing of radius $t = 1$ under the new $D$ metric.

**Proof.** To prove the theorem, it is sufficient to present an algorithm by which the constructed parity-check matrix will satisfy the conditions of Theorem 3.2. Let $\boldsymbol{h}_1, \boldsymbol{h}_2, \cdots, \boldsymbol{h}_n$ be the columns of the $\boldsymbol{H}$, where each $\boldsymbol{h}_i \in \mathbb{Z}_2^m$.

**Input data:** The dimension of the space: $n$.

**Output data:** An $m \times n$ parity-check matrix $\boldsymbol{H}$ that satisfies the conditions of Theorem 3.2.

**Algorithm steps:**
**Step 1: Check.** If the equality $n = 2^m - 1$ holds (where $m \geq 2$), then return the Hamming parity-check matrix, the construction of which is described 2nd chapter. Otherwise, go to Step 2.
**Step 2: Determination of parameters.** Find the minimum natural number $m$ (the number of parity-check bits) for which the inequality $2^m - 1 > n$ is satisfied. From the equality of the necessary condition for the sphere packing, $n_A \cdot n_B = 2^m - n - 1$, choose the following fixed values of the parameters:

$$n_A = 1, n_B = 2^m - n - 1.$$

From which the number of columns of the set $C$ is obtained: $n_C = n - n_A - n_B$. By substituting the values of the parameters $n_A$ and $n_B$, we get` $n_C = 2n - 2^m$.

**Step 3: Initialization.** Let $V = \mathbb{Z}_2^m \setminus \{\mathbf{0}\}$ be the set of all possible non-zero available vectors, the cardinality of which is $|V| = 2^m - 1$. Define empty sets of columns of the parity-check matrix for further filling: $\boldsymbol{H}_A := \emptyset$, $\boldsymbol{H}_B := \emptyset$ and $\boldsymbol{H}_C := \emptyset$.
**Step 4: Filling of subset $\boldsymbol{A}$.** From the set of vectors $V$, choose an arbitrary vector $\boldsymbol{h}_a$, assign it to the set $\boldsymbol{H}_A$, and remove it from the set $V$:

$$\boldsymbol{H}_A := \{\boldsymbol{h}_a\}, V := V \setminus \{\boldsymbol{h}_a\}.$$

**Step 5: Filling of subset $\boldsymbol{B}$.** Partition the remaining set $V$ into cosets according to the subgroup generated by the single element of the set $\boldsymbol{H}_A$, $\langle \boldsymbol{h}_a \rangle$ (each coset will have the form $\{\boldsymbol{\beta}, \boldsymbol{\beta} \oplus \boldsymbol{h}_a\}$). From these cosets, choose $n_B$ different cosets, and from each choose one vector, adding them to the set $\boldsymbol{H}_B$. Perform the selection as follows:

- **If** $n_C > 0$: choose $n_B$ arbitrary cosets and select one vector from each.
- **If** $n_C = 0$: In this case $n_B \geq 3$, it is necessary to choose two arbitrary cosets and take one vector from each (for example, $\boldsymbol{b}_1$ and $\boldsymbol{b}_2$). Then, from the coset corresponding to the sum of these vectors, the vector $\boldsymbol{b}_1 \oplus \boldsymbol{b}_2$ must be chosen. From the remaining cosets, choose one arbitrary vector each.

To ensure the 3rd condition of Theorem 3.2 (exclusion of forbidden sums), completely remove the chosen cosets from the available set $V$.
**Step 6. Filling of subset $\boldsymbol{C}$:** According to the equality $n_A \cdot n_B = 2^m - n - 1$, after removing the forbidden sums, exactly $n_C$ free vectors will remain in the set $V$. Assign all these vectors to the subset $C$: $\boldsymbol{H}_C := V$.
**Step 7. Matrix construction:** Combine the obtained sets, arranging them as columns of a single general matrix according to the corresponding positions: $\boldsymbol{H} = [\boldsymbol{H}_A \mid \boldsymbol{H}_C \mid \boldsymbol{H}_B]$. Return the matrix $\boldsymbol{H}$.

To ensure the correctness of the algorithm and complete the proof of Theorem 4.1, it is necessary to justify that the constructed matrix $\boldsymbol{H}$ satisfies the conditions of Theorem 3.2.

To verify the feasibility of Step 5, it is required to show that there is a sufficient number of cosets in the set of available vectors to form subsets $B$ and $C$. According to Step 5 of the algorithm, the space is partitioned into $2^{m-1}$ cosets with respect to the subgroup $\langle \boldsymbol{h}_a \rangle$, one of which is the subgroup itself. A total of $2^{m-1} - 1$ cosets remain available for selection. To construct subset $B$, it is required to use $n_B = 2^m - n - 1$ distinct cosets (taking one vector from each). This selection is possible because the required quantity does not exceed the number of available cosets:

$$n_B = 2^m - n - 1 \leq 2^m - 2^{m-1} - 1 = 2^{m-1} - 1,$$

which proves that we always have a sufficient number of cosets to select $n_B$ vectors. After selecting $n_B$ cosets, the number of unused (complete) cosets is:

$$2^{m-1} - 1 - n_B = 2^{m-1} - 1 - (2^m - n - 1) = n - 2^{m-1}.$$

Now let us consider subset $C$ constructed in Step 5, for which the required number of elements is $n_C = 2n - 2^m$. Since $C$ is formed by complete cosets (each containing 2 vectors), the number of cosets required for it will be half of the necessary elements:

$$\frac{n_C}{2} = \frac{2n - 2^m}{2} = n - 2^{m-1}.$$

The number of cosets required for subset $C$ exactly matches the number of remaining complete cosets in the system. Therefore, the execution of Step 5 is guaranteed.

At the beginning of the algorithm, the set of available non-zero vectors $V = \mathbb{Z}_2^m \setminus \{\mathbf{0}\}$ is formed. At each step, when choosing columns for sets $\boldsymbol{H}_A$, $\boldsymbol{H}_B$ and $\boldsymbol{H}_C$, these columns and their forbidden sums are immediately removed from the available set $V$. This guarantees that all columns of the final matrix $\boldsymbol{H}$ are pairwise distinct and non-zero, satisfying Condition 1.

The algorithm ensures the presence of linearly dependent columns satisfying Condition 2 by considering two possible cases. Case one: when $n_C > 0$. Since $n_C = 2n - 2^m$ is an even number, it is at least equal to 2. According to Step 5 of the algorithm, at least one complete coset will appear in the set $\boldsymbol{H}_C$. The two vectors of that coset, together with the vector $\boldsymbol{h}_a$, form a triplet whose sum is zero: $\boldsymbol{h}_a \oplus \boldsymbol{\beta} \oplus (\boldsymbol{\beta} \oplus \boldsymbol{h}_a) = \mathbf{0}$.
Case two: when $n_C = 0$. In this case, from the equation $n_C = 2n - 2^m = 0$, we get that $n = 2^{m-1}$. Substituting the value of n into the formula fixed in step 2, we have:

$$n_B = 2^m - n - 1 = 2^m - 2^{m-1} - 1 = 2^{m-1} - 1.$$

Since by the problem statement $m \geq 3$, then $n_B \geq 2^2 - 1 = 3$. Having at least 3 cosets, due to the construction of the 5th step of the algorithm, we choose for the set $\boldsymbol{H}_B$ such vectors $\boldsymbol{b}_1$, $\boldsymbol{b}_2$ and $\boldsymbol{b}_3 = \boldsymbol{b}_1 \oplus \boldsymbol{b}_2$ that satisfy the condition $\boldsymbol{b}_1 \oplus \boldsymbol{b}_2 \oplus \boldsymbol{b}_3 = \mathbf{0}$ within the set $\boldsymbol{H}_B$ itself.

When the algorithm selects a vector $\boldsymbol{h}_b$ for the set $\boldsymbol{H}_B$, the entire class corresponding to it is removed from the set of available vectors. Therefore, the forbidden sum $\boldsymbol{h}_b \oplus \boldsymbol{h}_a$ vector can no longer appear in the verification matrix, satisfying condition 3. Since the set $\boldsymbol{H}_A$ has only one element, therefore condition 4 is automatically satisfied.

Thus, the verification matrix constructed by the algorithm satisfies all the conditions of Theorem 3.2, which proves Theorem 4.1. ■

**Example 4.1** (Constructing a perfect sphere packing accordint to Theorem 4.1)**.** To demonstrate the algorithm for constructing a matrix, let's start by choosing the length of the code. Suppose we want to construct a code whose length is $n = 13$.
**Step 1: Check.** Since the number $n = 13$ cannot represented in the form $2^m - 1$, let's move on to step 2.
**Step 2: Determination of parameters.** Since it is required to find the smallest natural number $m$ that will satisfy the inequality $2^{m-1} \leq n < 2^m - 1$, in the case of $n = 13$ we get the value $m = 4$. Now let's calculate the powers of the primes according to the formulas of our algorithm:

$$n_A = 1,$$

$$n_C = 2n - 2^m = 10,$$

$$n_B = 2^m - n - 1 = 2.$$

**Step 3: Initialization.** We define the set of vectors as the set of all nonzero vectors of the $m$-dimensional space, $V = \mathbb{Z}_2^m \setminus \{\mathbf{0}\}$, which has $2^4 - 1 = 15$ elements. Iinitially, our target sets are empty: $\boldsymbol{H}_A \coloneqq \emptyset$, $\boldsymbol{H}_B \coloneqq \emptyset$ and $\boldsymbol{H}_C \coloneqq \emptyset$.

**Step 4: Filling of subset $\boldsymbol{A}$.** We choose any non-zero element from the space $V$, for example,

$$\boldsymbol{h}_a = (0,0,0,1).$$

**Step 5: Filling of subset $\boldsymbol{B}$.** In the set $V$ the remaining 14 elements are divided into $(16-2)/2 = 7$ cosets by the subgroup $\langle \boldsymbol{h}_a \rangle$:

$$\langle \boldsymbol{h}_a \rangle \oplus (0,0,1,0) = \{(0,0,1,0),(0,0,1,1)\},$$

$$\langle \boldsymbol{h}_a \rangle \oplus (0,1,0,0) = \{(0,1,0,0),(0,1,0,1)\},$$

$$\langle \boldsymbol{h}_a \rangle \oplus (0,1,1,0) = \{(0,1,1,0),(0,1,1,1)\},$$

$$\langle \boldsymbol{h}_a \rangle \oplus (1,0,0,0) = \{(1,0,0,0),(1,0,0,1)\},$$

$$\langle \boldsymbol{h}_a \rangle \oplus (1,0,1,0) = \{(1,0,1,0),(1,0,1,1)\},$$

$$\langle \boldsymbol{h}_a \rangle \oplus (1,1,0,0) = \{(1,1,0,0),(1,1,0,1)\},$$

$$\langle \boldsymbol{h}_a \rangle \oplus (1,1,1,0) = \{(1,1,1,0),(1,1,1,1)\}.$$

Since $n_C > 0$, we need to choose $n_B = 2$ arbitrary vectors from different cosets. Let us select the representatives of the first and second cosets:

$$\boldsymbol{h}_{b_1} = (0,0,1,0), \boldsymbol{h}_{b_2} = (0,1,0,0), \boldsymbol{H}_B = \{(0,0,1,0),(0,1,0,0)\}.$$

The cosets $\langle \boldsymbol{h}_a \rangle \oplus (0,0,1,0)$ and $\langle \boldsymbol{h}_a \rangle \oplus (0,1,0,0)$ are completely removed from the set of available vectors.

**Step 6: Filling of subset $\boldsymbol{C}$.** We need to take $n_C = 10$ vectors. Since exactly 5 cosets remain, each containing 2 elements, we simply take all the remaining 10 vectors. Note that the set $\boldsymbol{H}_C$ contains complete cosets, each of which, together with the vector $\boldsymbol{h}_a$, is linearly dependent. This ensures the presence of 3 linearly dependent columns in the subset of indices $I_1$, thereby satisfying the second condition of Theorem 3.2.

**Step 7: Matrix construction.** By combining the selected sets of columns in sequence, we obtain the $m \times n$ parity-check matrix $\boldsymbol{H}$:

$$\boldsymbol{H} = \begin{pmatrix} 0 & 0 & 0 & 1 & 1 & 1 & 1 & 1 & 1 & 1 & 1 & 0 & 0 \\ 0 & 1 & 1 & 0 & 0 & 0 & 0 & 1 & 1 & 1 & 1 & 0 & 1 \\ 0 & 1 & 1 & 0 & 0 & 1 & 1 & 0 & 0 & 1 & 1 & 1 & 0 \\ 1 & 0 & 1 & 0 & 1 & 0 & 1 & 0 & 1 & 0 & 1 & 0 & 0 \end{pmatrix}.$$

The construction algorithm of Theorem 4.1 provides a clear path for constructing a parity-check matrix based on fixed values of $n_A$ and $n_B$. However, there are also other parameters for which a parity-check matrix satisfying the conditions of Theorem 3.2 can be constructed. In particular, for a code length of $n = 13$, it is possible to consider other pairs of parameters, for example, $(n_A, n_B) = (2,9)$ or $(3,6)$. which require $m = 5$ check bits to ensure the partition conditions. Let us consider another example for the parameters $n_A = 3$, $n_B = 6$ and $m = 5$.

**Example 4.2** (Alternative parity-check matrix construction)**.** Suppose we want to construct a code with a length of 13, given the parameters $n_A = 3$, $n_B = 6$ and $m = 5$. To ensure the required partition, we select the set $\boldsymbol{H}_A$ such that, together with the zero vector, it forms a 2-dimensional subspace. Choosing a smaller dimension is not possible, and for a larger dimension, the number of cosets will not be sufficient to select $n_B$ columns. Let us consider the following subspace:

$$U = \{(0,0,0,0,0), (0,0,0,0,1), (0,0,0,1,0), (0,0,0,1,1)\}.$$

As the vectors of the set $\boldsymbol{H}_A$, we take the non-zero vectors of the subspace $U$. The entire space $\mathbb{Z}_2^5$ is partitioned with respect to the subspace $U$ into 8 cosets, each containing exactly 4 vectors. To exclude forbidden sums, we must choose one representative from different cosets as elements of the set $\boldsymbol{H}_B$. After selection, these cosets are completely removed from the set of available vectors. Let us select the $n_B = 6$ vectors in the following way:

$$\boldsymbol{H}_B = \{(0,0,1,0,0), (0,1,0,0,0), (0,1,1,0,0), (1,0,0,0,0), (1,0,1,0,0), (1,1,0,0,0)\}.$$

Such a choice also guarantees the fourth condition of Theorem 3.2, ensuring that forbidden sums of the type $\boldsymbol{h}_{a_1} \oplus \boldsymbol{h}_{b_1} = \boldsymbol{h}_{a_2} \oplus \boldsymbol{h}_{b_2}$ cannot occur. Suppose the contrary: that such a choice of columns for the set $\boldsymbol{H}_B$ produces forbidden sums that violate the fourth condition of Theorem 3.2. By rearranging the terms as $\boldsymbol{h}_{a_1} \oplus \boldsymbol{h}_{a_2} = \boldsymbol{h}_{b_1} \oplus \boldsymbol{h}_{b_2}$, we arrive at a contradiction. The left-hand side of the equation necessarily belongs to the subspace $U$, whereas the right-hand side can never belong to $U$ because the elements of the set $\boldsymbol{H}_B$ are chosen from strictly distinct cosets. Thus, we have exhausted 6 cosets, leaving only one remaining available coset, which contains exactly $n_C = 4$ vectors. We select all the vectors of this complete coset to form the set $\boldsymbol{H}_C$:

$$\boldsymbol{H}_C = \{(1,1,1,0,0), (1,1,1,0,1), (1,1,1,1,0), (1,1,1,1,1)\}.$$

By combining the selected sets of columns in sequence, we obtain the parity-check matrix $\boldsymbol{H}$ that satisfies the conditions of Theorem 3.2:

$$\boldsymbol{H} = \begin{pmatrix} 0 & 0 & 0 & 1 & 1 & 1 & 1 & 0 & 0 & 0 & 1 & 1 & 1 \\ 0 & 0 & 0 & 1 & 1 & 1 & 1 & 0 & 1 & 1 & 0 & 0 & 1 \\ 0 & 0 & 0 & 1 & 1 & 1 & 1 & 1 & 0 & 1 & 0 & 1 & 0 \\ 0 & 1 & 1 & 0 & 0 & 1 & 1 & 0 & 0 & 0 & 0 & 0 & 0 \\ 1 & 0 & 1 & 0 & 1 & 0 & 1 & 0 & 0 & 0 & 0 & 0 & 0 \end{pmatrix}.$$

As previously mentioned, for the case where $n = 13$ and $m = 5$, there exists another pair of parameters, namely $(n_A, n_B) = (2,9)$. However, unlike the choice of parameters in Example 4.2, it is impossible to construct a parity-check matrix satisfying the conditions of Theorem 3.2 using these specific parameters. Since $n_A = 2$, the set $\boldsymbol{H}_A$ consists of two distinct non-zero vectors, $\boldsymbol{H}_A = \{\boldsymbol{h}_{a_1}, \boldsymbol{h}_{a_2}\}$. Let us consider the subspace spanned by these vectors:

$$U = \langle \boldsymbol{h}_{a_1}, \boldsymbol{h}_{a_2} \rangle = \{\boldsymbol{0}, \boldsymbol{h}_{a_1}, \boldsymbol{h}_{a_2}, \boldsymbol{h}_{a_1} \oplus \boldsymbol{h}_{a_2}\}.$$

Since the vectors are linearly independent, the dimension of this subspace is 2, and it contains exactly 4 elements. The entire space $\mathbb{Z}_2^5$ is partitioned with respect to the subspace $U$ into $32/4 = 8$ cosets. To exclude forbidden sums, the elements of the set $\boldsymbol{H}_B$ cannot belong to the subspace $U$. Therefore, the $n_B = 9$ elements of $\boldsymbol{H}_B$ must be distributed across the remaining 7 available cosets. According to the Pigeonhole Principle (Dirichlet's box principle), if 9 elements are distributed among 7 cosets, then at least one coset will contain two elements from $\boldsymbol{H}_B$. Let these be the vectors $\boldsymbol{h}_{b_1}$ and $\boldsymbol{h}_{b_2}$ (where $\boldsymbol{h}_{b_1} \neq \boldsymbol{h}_{b_2}$). Since the vectors $\boldsymbol{h}_{b_1}$ and $\boldsymbol{h}_{b_2}$ lie within the same coset, their sum must belong to the subspace $U$:

$$\boldsymbol{h}_{b_1} \oplus \boldsymbol{h}_{b_2} \in U.$$

Since they are distinct elements, their sum is non-zero. Therefore, their sum must be equal to one of the remaining three non-zero vectors of $U$:

$$\boldsymbol{h}_{b_1} \oplus \boldsymbol{h}_{b_2} = \boldsymbol{h}_{a_1},$$
$$\boldsymbol{h}_{b_1} \oplus \boldsymbol{h}_{b_2} = \boldsymbol{h}_{a_2},$$
$$\boldsymbol{h}_{b_1} \oplus \boldsymbol{h}_{b_2} = \boldsymbol{h}_{a_1} \oplus \boldsymbol{h}_{a_2}.$$

All three cases immediately lead to forbidden sums, which violate the primary conditions of the matrix construction (Theorem 3.2).

## 5. Conclusion

By introducing the subsets of indices ($A$, $B$ and $C$), the cardinality of the ball ceases to depend solely on the total dimension $n$ and becomes a parameterizable quantity (depending on the parameters $n_A$ and $n_B$). Such a distance metric makes it possible to search for balls whose cardinalities satisfy the necessary condition for partitioning in a space of any dimension. Based on this property, the necessary and sufficient conditions for the partitioning of a Boolean space have been derived, where the latter is expressed through the linear dependencies and independencies of the columns of the parity-check matrix. As the most significant result of this work, it has been mathematically proven that a partition exists for a Boolean space of any dimension $n$.

Thus, the proposed generalized distance expands the possibilities for perfect sphere packing, rendering possible what was unattainable within the classical Hamming metric.

*This research did not receive any specific grant from funding agencies in the public, commercial, or not-for-profit sectors.*

## 6. Declarations

**Funding and Conflict of Interest:** The authors declare that no funds, grants, or other support were received during the preparation of this manuscript, and they have no relevant financial or non-financial interests to disclose.

**Author Contributions:** All authors contributed to the study conception and design. Material preparation and analysis were performed by T. Soghomonyan and Zh. Margaryan. The first draft of the manuscript was written by T. Soghomonyan, and all authors reviewed and approved the final version.

**Data Availability Statement:** Data sharing is not applicable to this article as no datasets were generated or analyzed during the current study. All mathematical results and algorithms are fully presented within the manuscript.

## 7. References

[1] Al-Bassam, S., Venkatesan, S., Al-Muhammadi, S.: New single asymmetric error-correcting codes. IEEE Trans. Inf. Theory 43(5), 1619–1623 (1997). https://doi.org/10.1109/18.623162

[2] AlBdaiwi, B., Horak, P., Milazzo, L.: Enumerating and decoding perfect linear Lee codes. Des. Codes Cryptogr. 52(2), 155–162 (2009). https://doi.org/10.1007/s10623-009-9273-3

[3] Astola, J.: On perfect codes in the Lee metric. Ann. Univ. Turku Ser. A I 176, 56 (1978)

[4] Berlekamp, E.R.: Algebraic Coding Theory. McGraw-Hill, New York (1968)

[5] Brouwer, A.E., Shearer, J.B., Sloane, N.J.A., Smith, W.D.: A new table of constant weight codes. IEEE Trans. Inf. Theory 36(6), 1334–1380 (1990). https://doi.org/10.1109/18.59932

[6] Brualdi, R.A., Graves, J.S., Lawrence, K.M.: Codes with a poset metric. Discrete Math. 147(1-3), 57–72 (1995). https://doi.org/10.1016/0012-365X(94)00228-B

[7] Chen, C.L., Hsiao, M.Y.: Error-correcting codes for semiconductor memory applications: A state-of-the-art review. IBM J. Res. Dev. 28(2), 124–134 (1984). https://doi.org/10.1147/rd.282.0124

[8] Chiang, J.C.Y., Wolf, J.K.: On channels and codes for the Lee metric. Inf. Control 19(2), 159–173 (1971). https://doi.org/10.1016/S0019-9958(71)90791-1

[9] Constantin, S., Rao, T.R.N.: On the theory of binary asymmetric error-correcting codes. Inf. Control 40(1), 20–36 (1979). https://doi.org/10.1016/j.ffa.2018.10.004

[10] Dougherty, S.T., Kim, J.-L., Solé, P.: Open problems in coding theory. Contemp. Math. 634, 79–99 (2015). https://doi.org/10.1090/conm/634/12692

[11] Ghandilyan, A.Kh., Margaryan, Z.G.: One error and in the given interval two error correcting code for the additive communication channel. Math. Probl. Comput. Sci. 32, 14–22 (2009)

[12] Golomb, S.W., Welch, L.R.: Algebraic coding and the Lee metric. In: Mann, H.B. (ed.) Error Correcting Codes, pp. 175–189. Wiley, New York (1968)

[13] Guruswami, V., Rudra, A., Sudan, M.: Essential Coding Theory. University at Buffalo, Buffalo (2022)

[14] Hamming, R.W.: Error detecting and error correcting codes. Bell Syst. Tech. J. 29(2), 147–160 (1950). https://doi.org/10.1002/j.1538-7305.1950.tb00463.x

[15] Hill, R.: A First Course in Coding Theory. Oxford University Press, Oxford (1986)

[16] Jain, S., Nam, K.-B., Lee, K.-S.: On some perfect codes with respect to Lee metric. Linear Algebra Appl. 405, 104–120 (2005). https://doi.org/10.1016/j.laa.2005.03.001

[17] Krotov, D.S., Sotnikova, E.V.: Embedding in q-ary 1-perfect codes and partitions. Discrete Math. 338(11), 1856–1859 (2015). https://doi.org/10.1016/j.disc.2015.04.014

[18] Lee, C.Y.: Some properties of non-binary error correcting codes. IEEE Trans. Inf. Theory 4(2), 77–82 (1958). https://doi.org/10.1109/TIT.1958.1057446

[19] Lovisolo, L.M., da Silva, E.A.B.: Generalized Hamming distance. IEE Proc. Vis. Image Signal Process. 148(1), 39–46 (2001). https://doi.org/10.1023/A:1020499411651

[20] MacWilliams, F.J., Sloane, N.J.A.: The Theory of Error-Correcting Codes. North-Holland Publishing Company, Amsterdam (1977)

[21] Mahdaviani, K., Shahidi, S., Haddadi, S., Ardakani, M., Tellambura, C.: Improving the sphere-packing bound for binary codes over memoryless symmetric channels. In: Proceedings of the 47th Annual Allerton Conference on Communication, Control, and Computing (Allerton), pp. 1041–1048 (2009). https://doi.org/10.1109/ALLERTON.2009.5394906

[22] Niederreiter, H.: Point sets and sequences with small discrepancy. Monatsh. Math. 104(4), 273–337 (1987). https://doi.org/10.1007/BF01294651

[23] Niederreiter, H.: A combinatorial problem for vector spaces over finite fields. Discrete Math. 96(3), 221–228 (1991). https://doi.org/10.1016/0012-365X(91)90315-S

[24] Szabó, S.: On mosaics consisting of multidimensional crosses. Acta Math. Acad. Sci. Hung. 38(1-4), 191–203 (1981). https://doi.org/10.1007/BF01917533
[25] Tietäväinen, A.: On the nonexistence of perfect codes over finite fields. SIAM J. Appl. Math. 24(1), 88–96 (1973). https://doi.org/10.1137/0124010
[26] Varshamov, R.R.: A class of codes for asymmetric channels and a problem from the additive theory of numbers. Soviet Math. Dokl. 14, 127–131 (1973). https://doi.org/10.1109/TIT.1973.1054954